\documentstyle[amsfonts,twoside,12pt]{article}
\input{psfig}
\def\SBIMSMark#1#2#3{
 \font\SBF=cmss10 at 10 true pt
 \font\SBI=cmssi10 at 10 true pt
 \setbox0=\hbox{\SBF Stony Brook IMS Preprint \##1}
 \setbox2=\hbox to \wd0{\hfil \SBI #2}
 \setbox4=\hbox to \wd0{\hfil \SBI #3}
 \setbox6=\hbox to \wd0{\hss
             \vbox{\hsize=\wd0 \parskip=0pt \baselineskip=10 true pt
                   \copy0 \break%
                   \copy2 \break%
                   \copy4 \break}}
 \dimen0=\ht6   \advance\dimen0 by \vsize \advance\dimen0 by 8 true pt
                \advance\dimen0 by -\pagetotal
 \dimen2=\hsize \advance\dimen2 by .25 true in
  \openin2=publishd.tex
  \ifeof2\setbox0=\hbox to 0pt{}
  \else 
     \setbox0=\hbox to 3.1 true in{
                \vbox to \ht6{\hsize=3 true in \parskip=0pt  \noindent  
                \input publishd.tex 
                \vfill}}
  \fi
  \closein2
  \ht0=0pt \dp0=0pt
 \ht6=0pt \dp6=0pt
 \setbox8=\vbox to \dimen0{\vfill \hbox to \dimen2{\copy0 \hss \copy6}}
 \ht8=0pt \dp8=0pt \wd8=0pt
 \copy8
 \message{*** Stony Brook IMS Preprint #1, #2 ***}
}

\def  \R{{\hbox{$\Bbb R$}}}
\def  \d{{\Bbb D}}

\def \a{{\overline \d}}
\def \aa{{\overline \d}^2}
\def \dn{{\aa}\times\dots\times {\aa}}
\def \NI{\noindent}
\def \H{{\rm Homeo}({\aa}, \partial {\aa})}

\def \HU{{\rm Homeo}({\aa}, \partial {\aa} ,\lambda_1, \dots,
\lambda_n)}
\def \K{{\rm Homeo}({\aa}, \partial {\aa} ,\mu_1, \dots, \mu_n)}

\def \D{{\rm Diff}^1({\aa}, \partial {\aa} )}

\newtheorem{thm}{Theorem}
\newtheorem{lem}{Lemma}
\newtheorem{prop}{Proposition}
\newtheorem{cor}{Corollary}
\newenvironment{pf}{\noindent {\bf Proof:}}{\nolinebreak \hfill $\Box$}

\title{ Dynamical Cocycles with Values \\ in the Artin Braid Group}
\author{J.-M. Gambaudo and E. E.  P\'ecou}
\date{}

\begin{document}
\thispagestyle{empty}

\maketitle
\SBIMSMark{1997/5}{April 1997}{}
\NI{\bf Abstract:}
By considering the way an $n$-tuple of points in the 2-disk are linked
together under iteration of an orientation preserving diffeomorphism,
we
construct a dynamical cocycle with values in the Artin braid group. We
study
the asymptotic properties of this cocycle and derive a series of
topological
invariants for the diffeomorphism which enjoy rich properties.

\bigskip\bigskip\bigskip

\section{Introduction}

When one is concerned with the study of a group of automorphisms on a
probability space $(X,{\cal B}, \mu)$, the knowledge of some cocycle
associated to this group is useful. Indeed, cocycles carry into a
simple, well understood, target group $G$ the main dynamical features.
By using subadditive functions on $G$,  for instance left-invariant
metrics, one is able to define asymptotic invariants with relevant
dynamical properties.  As an illustration, see the Oseledec theory of
the Lyapounov exponents (\cite{Oseledec}).

\bigskip

 In this paper, we start by recalling  some basic definitions related
 to cocycles in general and we state some nice asymptotic properties
 when a subadditive function on the target group is given. Namely, we
 show the existence of asymptotic invariants and give conditions for
 their topological invariance (section $2$).

\NI Section $3$ deals with the study of orientation preserving
$C^1$-diffeomorphisms of the $2$-disk which preserve a non atomic
measure.  Given a $n$-tuple of distinct points in the disk, we
construct on the group of diffeomorphisms we are considering  a
cocycle  with values in the Artin braid group $B_n$. Indeed, we show
that given a diffeomorphism there is a well defined way to associate a
braid to a $n$-tuple of orbits. This construction generalizes a
well-known construction for periodic orbits.

\NI  Despite the relative complexity of the Artin braid group, there
are naturally defined subadditive functions on $B_n$. Therefore, we can
consider the asymptotic invariants associated to these cocycles. These
invariants turn to be topological invariants that we can relate to
other well known quantities such as the Calabi invariant and the
topological entropy (section 4). In this last case, we generalize a
result of minoration of the entropy by Bowen (see \cite{Bowen}).

\section{Cocycles}
\subsection {Asymptotic behaviour}
Let $(X,{\cal B}, \mu)$ be a probability space,
 $Aut (X, \mu)$, the group of all its automorphisms, i.e. invertible
 measure preserving transformations and  $G$ a topological group (with
 identity element $e$).
 For any subgroup $\Gamma$ of $Aut(X, \mu)$, we say that a  measurable
 map
$\alpha \,:\, X\times \Gamma \to G$
is a {\it cocycle} of the dynamical system $(X, {\cal B}, \mu, \Gamma)$
with values in $G$ if for all $\gamma_1$ and $\gamma_2$ in $\Gamma$ and
for $\mu$-a.e. $x$ in $X$:
$$ \alpha(x, \gamma_1\gamma_2)\,=\, \alpha (x,
\gamma_1)\,\alpha(\gamma_1 x, \gamma_2).$$
\NI A continuous map $\theta : G\to \R^+$, which satisfies, for all
$g_1$ and $g_2$ in $G$:
$$\theta (g_1g_2)\,\,\leq \,\,\theta (g_1) \,+\, \theta(g_2)$$
is called a {\it subadditive function} on $G$.

\NI For example, if the group $G$ is equipped with a right (or left)
invariant metric $d$, then it is clear that the map:
$$\begin{array}{rcl}
	G & \to & \R^+\\
	x & \mapsto & d(g , e)
\end{array}$$
is a subadditive function.

\bigskip

The following Theorem is a straightforward application of the
subadditive ergodic theorem (see for instance \cite {Walter}).

\begin{thm}
Let  $\alpha$ be a cocycle of the dynamical system $(X, {\cal B}, \mu,
\Gamma)$ with values in a group $G$  and $\theta$ a subadditive
function on $G$.
 Assume that there exists an automorphism $\gamma$ in $\Gamma$ such
 that the map:
$$\begin{array}{rcl}
	X & \to & \R^+\\
	x & \mapsto & \theta(\alpha (x, \gamma))
\end{array}$$
is integrable. Then:

\begin{enumerate}
\item For $\mu$-a.e. $x$ in $X$, the sequence ${{1}\over{n}}
\theta(\alpha (x, \gamma^n))$ converges, when $n$ goes to $+\infty$ to
a limit denoted by  $\Theta(x, \gamma, \alpha)$.

\item The map:
$$\begin{array}{rcl}
   X & \to & \R^+\\
   x & \mapsto & \Theta(x, \gamma, \alpha)
\end{array}$$
 is integrable and we have :
$$\lim_{n\to +\infty} {{1}\over{n}}\int \theta(\alpha(x, \gamma^n))d\mu
\,\,=\,\,
 \inf_{n\to +\infty} {{1}\over{n}}\int \theta(\alpha(x,\gamma^n))d\mu
 \,\,=\,\,
\int \Theta(x, \gamma, \alpha) d\mu.$$
\end{enumerate}
\end{thm}

\NI  We denote  by $\Theta_\mu(\gamma, \alpha)$ the quantity $\int
\Theta(x, \gamma, \alpha) d\mu$.

\subsection {Invariance property}
Let $\alpha$ and $\beta$ be two cocycles of the dynamical system  $(X,
{\cal B}, \mu, \Gamma)$  with value in a group $G$. We say that
$\alpha$ and $\beta$  are {\it cohomologous} if there exists a
measurable function $\phi:\,  X \to G$ such that for all $\gamma$ in
$\Gamma$ and for $\mu$-a.e. $x$ in $X$:

$$ \beta(x, \gamma)\,=\, \phi(x)\,\alpha(x, \gamma)\,
(\phi(\gamma\,x))^{-1}.$$

\bigskip

Consider two probability spaces $(X_1, {\cal B}_1, \mu_1)$ and
 $(X_2, {\cal B}_2, \mu_2)$ and  an isomorphism $h: X_1 \to X_2$,
( i.e. an invertible transformation which carries the measure $\mu_1 $
onto the measure $\mu_2$). Let $\Gamma_1$ be a subgroup of $Aut(X_1,
\mu_1)$ and $\Gamma_2\,=\,h\Gamma_1 h^{-1}$. Given a cocycle
$\alpha_1$  on $(X_1, {\cal B}_1, \mu_1, \Gamma_1)$, we denote by
$h\circ\alpha_1$ the cocycle defined on $(X_2, {\cal B}_2, \mu_2,
\Gamma_2 )$ by:
$$
h\circ\alpha_1(x_2, \gamma_2)\,=\, \alpha_1(h^{-1}(x_2), h^{-1}\gamma_2
h),
$$
for all $\gamma_2 \in \Gamma_2$ and for all $x_2\in X_2$.

\bigskip

Given two dynamical systems $(X_1, {\cal B}_1, \mu_1, \Gamma_1)$ and
$(X_2, {\cal B}_2, \mu_2, \Gamma_2)$ and given two cocycles $\alpha_1$
and $\alpha_2$ defined respectively on the first and the second system
and with values in the same group $G$, we say that the two cocycles
$\alpha_1$ and $\alpha_2$ are {\it weakly equivalent} if there exists
an isomorphism   $h : X_1 \to X_2$ which satisfies :
\begin {enumerate}
\item $\Gamma_2 \,=\,   h\Gamma_1 h^{-1}$,
\item $\alpha_2$ and $h\circ \alpha_1$ are cohomologous.

\end{enumerate}

\bigskip

\begin{thm}
{Let  $\alpha_1$ and $\alpha_2$ be two cocycles defined on the
dynamical
systems $(X_1, {\cal B}_1, \mu_1, \Gamma_1)$ and $(X_2, {\cal B}_2,
\mu_2, \Gamma_2)$ respectively and with values in a group equipped with
a subadditive function  $\theta$
. Assume that $\alpha_1$ and $\alpha_2$ are  weakly equivalent through
an isomorphism $h$.

\NI If  $\gamma_1\in\,\Gamma_1$ and $\gamma_2 = h\gamma_1 h^{-1}$
satisfy that both maps:
$$x_1 \mapsto \theta\,(\alpha_1(x_1, \gamma_1))\,\,\,\,\,\,\,\,\,\,{\it
and}\,\,\,\,\,\,\,\,\,\,x_2 \mapsto \theta\,(\alpha_2(x_2, \gamma_2))$$
are integrable respectively on $X_1$ and $X_2$, then:
$${\Theta}_{\mu_1}(\gamma_1, \alpha_1)\,=\, {\Theta}_{\mu_2}(\gamma_2,
\alpha_2).$$
}\end{thm}
\NI Before proceeding to the proof of the theorem, let us recall a
basic lemma in ergodic theory:
\begin{lem} Let $(X, {\cal B}, \mu)$ be a probability space, $\gamma$
an automorphism on $X$,  $G$ a group, $\theta$ a subadditive function
on $G$ and  $\phi: X\to G$ a measurable map.
 Then,  for $\mu$-a.e.  $x$ in $X$:
$$
\liminf_{n\to +\infty} \theta(\phi(\gamma^n x)) \,<\,\,+\infty.$$
\end{lem}

\medskip

\begin{pf} Let $N$ be a positive integer and
$${\cal E}_N\,\,=\,\,\{x\in X\,\,\vert\,\, \theta(\phi(x))\leq N\}.$$
\NI Since the map $\phi$ is measurable the set ${\cal E}_N$ is
measurable and has a positive measure for $N$ big enough. From the
Poincar\'e recurrence theorem (see for instance  \cite{Walter}) we know
that, for $\mu$-a.e $x$ in ${\cal E}_N$,   the orbit $\gamma^n(x)$
visits ${\cal E}_N$ infinitely often and thus:
$$\liminf_{n\to +\infty} \theta(\phi(\gamma^n x)) \,\leq\,\,N.$$
\NI Since the union $\cup_{N\geq0} {\cal E}_N$ is a set of
$\mu$-measure 1 in $X$, the lemma is proved.
\end{pf}
\bigskip

\NI {\bf Proof of Theorem 2:} The weak equivalence of the cocycles
$\alpha_2$ and $\alpha_1$ leads to the existence of some measurable
function $\phi : X_2 \to G$  such that, for all $n \geq 0$, for all
$\gamma_2$ in $\Gamma_2$ and for $\mu_2$-a.e. $x_2$ in $X_2$:
$$
\alpha_2(x_2, \gamma_2^n) \,=\, \phi(x_2)\,h\circ\alpha_1(x_2,
\gamma_2^n)\,(\phi(\gamma_2^nx_2))^{-1}.$$

\NI The  subadditivity property of the map $\theta : G\to \R^+$   gives
us:
$$
\vert \theta(\alpha_2(x_2, \gamma_2^n))\,-\, \theta(h\circ\alpha_1(x_2,
\gamma_2^n))\vert  \,\,\leq\,\,
 \theta(\phi(x_2))\,+\, \theta(\phi(\gamma_2^n x_2)^{-1}).
$$
\NI From Lemma 1, it follows that for all $\gamma_2$ in $\Gamma_2$, and
for $\mu_2$ a.e. $x_2$ in $X_2$:
$$
\liminf_{n\to+\infty}\vert \theta(\alpha_2(x_2, \gamma_2^n))\,-\,
\theta(h\circ\alpha_1(x_2, \gamma_2^n) )\vert \,<\, +\infty.$$
\NI Thus:
$$
\liminf_{n\to+\infty}\vert {{1}\over{n}}\theta(\alpha_2(x_2,
\gamma_2^n))\,-\,{{1}\over{n}} \theta(h\circ\alpha_1(x_2, \gamma_2^n)
)\vert \,=\,0.
$$
Since the map $\theta\,(\alpha_i(\,.\,, \gamma_i))$ is integrable for
$i=1$ and $i=2$, it results from Theorem $1$ that for $\mu_i$-a.e. $x$
in $X_i$:
$$
\lim_{n\to +\infty}{{1}\over{n}} \theta\,(\alpha_i(x, \gamma_i^n))  =
\Theta (x, \gamma_i, \alpha_i)
$$
Recalling that the definition of $h\circ\alpha_1(x_2, \gamma_2^n)$ is
$\alpha_1(h^{-1}(x_2), h^{-1}\gamma_2^n h)$, we get that, for
$\mu_2$-a.e. $x_2$ in $X_2$:
$$
\Theta(h^{-1}(x_2), \gamma_1, \alpha_1)\,\,=\,\, \Theta(x_2, \gamma_2,
\alpha_2).$$
By integrating this equality with respect to $\mu_2$, we get:
$${\Theta}_{\mu_1}(\gamma_1, \alpha_1) \,\,=\,\,
{\Theta}_{\mu_2}(\gamma_2, \alpha_2).$$

\section{Braids and Dynamics}

\subsection{The Artin braid group}
For a given integer $n>0$, the Artin braid group $B_n$ is the group
defined by the generators $\sigma_1, \dots , \sigma_{n-1}$ and the
relations:
$$\begin{array}{rcll}
\sigma_i\sigma_j & = & \sigma_j\sigma_i, & \vert i-j\vert \geq
2,\,\,\,\,\, 1\leq i, j\leq n-1\\
\sigma_i\sigma_{i+1}\sigma_i & = & \sigma_{i+1}\sigma_i\sigma_{i+1}, &
1\leq i\leq n-2.
\end{array}
$$
 An element in this group is called a {\it braid}.

\medskip

\NI A geometrical representation of the Artin braid group is given by
the following construction.

\NI Let $\a^2$ denote the closed unit disk in $\R^2$ centered at the
origin and $\d^2$ its interior. Let $Q_n = \{q_1, \dots , q_n\}$ be a
set of $n$ distinct points in $\d^2$ equidistributed on a diameter of
$\d^2$. We denote by $\d_n$ the $n$-punctured disk ${\a^2} \setminus
Q_n$.

\NI We define a collection of $n$ arcs in the cylinder $ {\d^2} \times
[0, 1]$:
$$\Gamma_i = \{(\gamma_i(t), t) \,,\,t\in\,[0,1]\}, \,\,\,\,\,\,\,
1\leq i \leq n$$
joining the points in $Q_n \times \{ 0\}$ to the points in $Q_n \times
\{1\}$ and  such that $\gamma_i(t) \neq \gamma_j(t)$ for $i\neq j$.

\NI We call $\Gamma = \cup_{i=0,\dots , n}\Gamma_i$ a {\it geometrical
braid} and say that two geometrical braids are equivalent if there
exists a continuous deformation from one to the other through
geometrical braids. The set of equivalence classes is the Artin braid
group. The composition law
 is given by concatenation as shown in Figure 1 and a  generator
 $\sigma_i$ corresponds to the geometrical braid shown in Figure 2.

\begin{figure}[ht]
\centerline{\psfig{figure=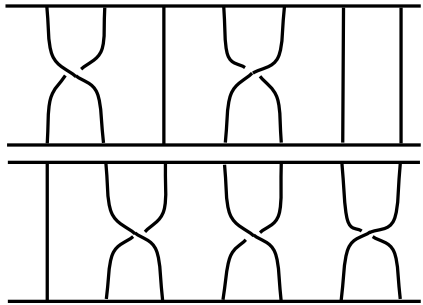,height=40mm}}
\caption{The concatenation of two geometrical braids}
\end{figure}

\begin{figure}[ht]
\centerline{\psfig{figure=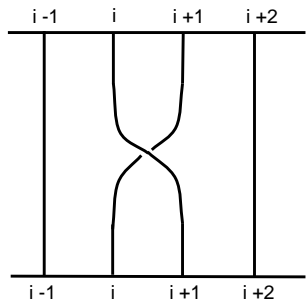,height=40mm}}
\caption{Representation of the $i$-th generator $\sigma_i$ of the Artin
braid group}
\end{figure}

\NI
Let $\H$ be the subset of homeomorphisms from $\a^2$ onto itself which
are the iudentity map in a neighborhood on the boundary of the disk.
The subgroup of elements  which let the set $Q_n$ globally invariant is
denoted by ${\rm Homeo}({\a^2}, \partial {\a^2}, Q_n)$.  The following
Theorem, due to J. Birman, shows the relation between braids and
dynamics.

\begin{thm}[\cite{Birman}]{The Artin braid group $B_n$ is isomorphic to
the  group of automorphisms of $\pi_1({\d_n})$ which are induced by
elements of ${\rm Homeo}({\a^2}, \partial {\a^2}, Q_n)$, that is to say
the group of isotopy classes of ${\rm Homeo}({\a^2}, \partial {\a^2},
Q_n)$.}
\end{thm}

\NI This isomorphism ${\cal I}$  can be described as follows. Let
$\{x_1, \dots , x_n\}$ be a basis for the free group $\pi_1({\d_n})$,
where $x_i$, for $1\leq i\leq n$ is represented by a simple loop which
encloses the boundary point $q_i$ (see Figure 3). A generator
$\sigma_i$ of $B_n$ induces on $\pi_1({\d_n})$ the following
automorphism ${\cal I}(\sigma_i)$ :

$${\cal I}(\sigma_i)\,\left \{ \begin{array}{lcl}
    x_i     & \mapsto & x_i\, x_{i+1}\, {x_i}^{-1}\\
    x_{i+1} & \mapsto & x_i\\
    x_j     & \mapsto & x_j \,\,\,\,\,\,{\rm if}\,\, \,\,j\,\neq\,
    i\,,\, i+1.
  \end{array} \right.$$

\NI The action of the Artin braid group on $\pi_1({\d_n})$ is a right
action. We denote by $wb$ the  image of  $w\in \pi_1({\d_n})$ under the
automorphism induced by the braid $b$.

\begin{figure}[ht]
\centerline{\psfig{figure=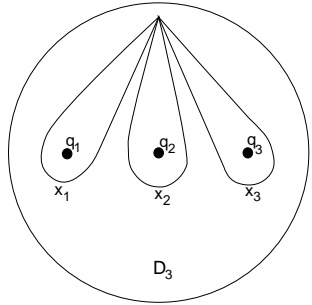,height=40mm}}
\caption{The generators of \protect{$\pi_1({\bf D}_n)$}}
\end{figure}

\bigskip

  We now introduce two standard subadditive functions on the Artin
  braid group. Given  a group presented  by a finite number of
  generators and relations and $g$ an element of the group,  we denote
  by $L(g)$ the minimal length of the word $g$ written with the
  generators  and their inverses.
We can define two subadditive functions on $B_n$ by  measuring lengths
of words either in the Artin braid group or in the fundamental group of
the punctured disk $\d_n$. More precisely, given  an  element $b$  in
$B_n$,  the first  subadditive function $\theta_1$ is defined by:
$$\theta_1(b) \,\,= \,\, L(b).$$
Notice that by setting $d(b, c)\,=\, L(bc^{-1})$ for all $b$ and $c$ in
$B_n$, we define a right invariant metric $d$ on $B_n$.

\NI
The second subadditive function $\theta_2$ is defined, for all $b$ in
$B_n$, by:
$$ \theta_2(b)\,\,=\,\, \sup_{i=1,\dots, n}\log L(x_ib),$$
where the $x_i$'s are the generators of the fundamental group of $\d_n$
defined as above.

\NI It is plain that the two subadditive functions  are related as
follows:
$$ \theta_2(b)\,\,\leq \,\,(\log 3)\theta_1(b),$$
for all $b$ in  $B_n$.

\bigskip
\NI In the sequel, we shall focus on a particular subgroup of the Artin
braid group, that we now define. Let $b$ be an element in $B_n$
represented by a geometrical braid:
$$\Gamma = \cup_{i=0,\dots,n}\{(\gamma_i(t), t) \subset {\d^2} \times
[0, 1], \,\,\, 1\leq i \leq n\}.$$
We say that $b$ is a {\it pure braid} if for all $i\,=\,1,\dots,n$,
$\gamma_i(0)\,=\,\gamma_i(1)$.

\NI We denote by ${\cal F}_n$ the set of all pure braids in $B_n$ (see
for instance \cite{Zieschang}) .

\medskip

\NI {\bf Remark:} The pure braid group can be interpreted as follows.
Consider a braid  $b$ in ${\cal F}_n$, and  $\Gamma = \cup_{i=0,\dots ,
n} \{(\gamma_i(t), t), t\in [0, 1]\}$  a geometrical braid whose
equivalence class is $b$. We can associate to $\Gamma$ the loop
$(\gamma_1(t), \dots ,\gamma_n(t))_{t\in [0, 1]}$ in the space ${{\bf
D^2}\times\dots\times {\bf D^2} }\setminus \Delta$, where $\Delta$
 is the {\it generalized diagonal }:
$$ \Delta \,\,=\,\, \{(p_1, \dots, p_n) \in  {{\bf
D^2}\times\dots\times {\bf D^2} },\,\vert \,\,\,\,\exists\,  i\neq j
\,\,\, {\rm and}\,\,\, p_i  =  p_j\}.
$$

\NI In other words, there exist an isomorphism $\cal J$ between the
pure braid group ${\cal F}_n$ and  the   fundamental group
 $\pi_1({{\bf D^2}\times\dots\times {\bf D^2} }\setminus \Delta)$.

\subsection{A cocycle with values in the  Artin braid group }

\NI Let $\phi$ be an element in $\H$, and let $P_n\,=\,(p_1,\dots,
p_n)$ be a $n$-tuple of pairwise disjoint  points  in $\d^2$. The
following procedure describes a natural way to associate a pure braid
in ${\cal F}_n$ to the datas $\phi$ and $P_n$ (see figure 4).
\begin{enumerate}
\item For $1\leq i\leq n$,  we   join the point $(q_i,  0)$ to the
point  $(p_i,{1\over3} )$ in the cylinder $ {\d^2} \times [0, {1 \over
3}]$ with  a segment.
\item For $1\leq i\leq n$, we   join the point $(p_i, {1\over3} )$ to
the point  $(\phi(p_i),{2\over3} )$ in the cylinder $ {\d^2} \times [{1
\over 3}, {2\over 3}]$ with the arc $(\phi_{(3t-1)}(p_i), t)_{t\in [{1
\over 3}, {2\over 3} ]}$, where $(\phi_\tau)_{\tau\in [0, 1]}$ is an
isotopy from the identity to $\phi$ in $\H$.
\item For $1\leq i\leq n$, we   join the point $(\phi(p_i), {2\over3}
)$ to the point  $(q_i, 1)$ in the cylinder $ {\d^2} \times [{2 \over
3}, 1]$ with  a segment.
\end{enumerate}
\NI The  concatenation of this sequence of arcs provides  a geometrical
braid. The equivalence class of this geometrical braid is  a braid in
${\cal F}_n$ that we denote by $\beta (P_n; \phi)$.

\begin{figure}[ht]
\centerline{\psfig{figure=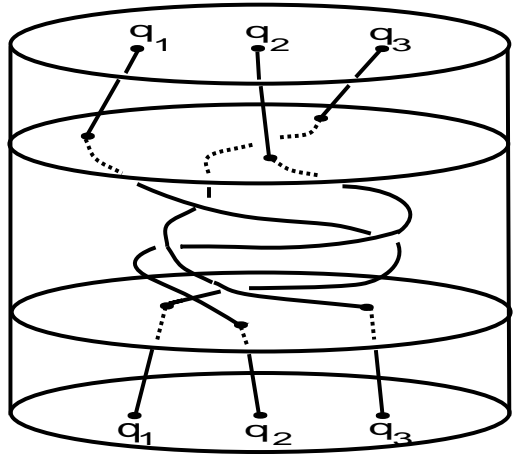,height=60mm}}
\caption{Construction of $\beta(P_n; \phi)$}
\end{figure}

\bigskip
\NI It is clear that the above procedure is well defined if and only if
for all $1\leq i<j\leq n$:
\begin{enumerate}
\item[{(i)}] the segment joining the points $(q_i, 0)$ and $(p_i,
{1\over3})$ does not intersect the segment joining the points $(q_j,
0)$ and $(p_j,{1\over3} )$, and

\item[{ (ii)}] the segment joining the points $(\phi(p_i),{2\over3} )$
and $(q_i, 1)$ does not intersect the segment joining the points
$(\phi(p_j),{2\over3} )$ and $(q_j, 1)$.
\end{enumerate}

\NI For $1\leq i\neq j\leq n$, consider in $\R^2\times\dots\times \R^2$
the codimension 2 plane $ P_{i,j}$  of points $(z_1, \dots, z_n)$ such
that $z_i=z_j$. Let $H_{i,j}$ be the hyperplane which contains $P_{i,
j}$  and the point $(q_1, \dots , q_n)$. The plane $P_{i, j}$ splits
the hyperplane $H_{i, j}$ in 2 components; we denote by $H^0_{i, j}$
the closure of the  component which does not contain the point $(q_1,
\dots,q_n)$.  Let $\Omega^{2n}$ be the open dense subset of $\dn$
defined by:
$$\Omega^{2n}\,\,=\,\, {\dn} \setminus \bigcup_{1\leq i , j \leq n}
H^0_{i, j}\cap {\dn}.$$

\NI We can reformulate  conditions (i) and (ii) as follows :  the braid
$\beta (P_n; \phi)$ is  defined if and only if both $P_n$ and
$(\phi,\dots,\phi) (P_n)$ belong to $\Omega^{2n}$.

\bigskip
\NI {\bf Remark 1:} This procedure is locally constant where it is
defined. More precisely \-  if $\beta (P_n; \phi)$ is defined and if
$P'_n$ is  close enough to $P_n$ in ${\dn}\setminus \Delta$ , then
  $\beta (P'_n; \phi)$ is also defined and  $\beta (P'_n; \phi)$ is
  equal to
$\beta (P_n; \phi)$.

\bigskip
\NI {\bf Remark 2:} Since the set $\H$ is contractible, this procedure
does not depend on the isotopy $\phi_t$.

\bigskip

\NI In order to construct a cocycle with values in the Artin braid
group, we also need to define the invariant measures we are going to
consider. We say that a $n$-tuple of probability measures
 $\lambda_1 , \dots , \lambda_n$ on $\a^2$ is {\it coherent} if the
 subset $\Omega^{2n}
$   has measure 1 with respect to the product measure
$\lambda_1\times\dots\times\lambda_n$.  Notice that this is the case
when the measures $\lambda_i$ are non atomic.

\begin{lem} Let  $\phi$ be a map in $\H$ which preserves the coherent
probability  measures $\lambda_1,\dots , \lambda_n$. Then,  the map :

$$
\begin{array}{rll}
		  {\dn} & \longrightarrow & {\cal F}_n\\
		    P_n & \longmapsto     & \beta (P_n; \phi)
\end{array}
$$
is measurable with respect to the product measure $\lambda_1 \times
\dots \times \lambda_n$.
\end{lem}

\begin{pf} The map is continuous on $\Omega^{2n}$ which has measure 1
with respect to the product measure
$\lambda_1\times\dots\times\lambda_n$.
\end{pf}

\bigskip

There is a natural injection $j$ from the set of maps from $\a^2$ into
itself into the set of maps from $\dn$ into itself. Namely:
$$j(\phi)\,=\,(\phi,\dots,\phi).$$
\NI Let $\lambda_1,\dots ,\lambda_n$ be coherent probability measures
and $\HU$ the set of maps in $\H$ which preserve the measures
$\lambda_1,\dots ,\lambda_n$.

\NI We call $\Gamma$ the image of  $\HU$ by $j$ in $Aut({\dn},
\lambda_1\times\dots\times\lambda_n)$.

\NI The following proposition is straightforward:

\begin{prop} The map:
$$
\begin{array}{rll}
		{\dn}\times \Gamma & \longrightarrow &  {\cal F}_n\\
		    (P_n,j(\phi) ) & \longmapsto     & \beta(P_n; \phi)
\end{array}
$$
is a cocycle of the dynamical system $({\dn}, {\cal B}({\dn}),
\lambda_1\times \dots\times \lambda_n, \Gamma)$ with values in the
group ${\cal F}_n$ (here ${\cal B}(\dn)$ stands for the  Borel
$\sigma$-algebra of $\dn$).
\end{prop}

\subsection{Asymptotic limits and invariance}
Let $\D$ denote the set of $C^1$-diffeomorphisms from $\a^2$ into
itself which are the identity in a neighborhood of the boundary. The
following result is fundamental for our purpose.

\begin{prop} If the map $\phi$ is in $\D$ then there exists a positive
number $K(\phi, n)$ such that for all  $P_n$ where   $\beta (P_n;
\phi)$ is defined :
$$ \theta_1(\beta(P_n;  \phi)) \leq K(\phi, n).$$
\end{prop}

\begin{pf} From the isomorphism ${\cal J}$ between ${\cal F}_n$ and
$\pi_1({{\bf D^2}\times\dots\times {\bf D^2} }\setminus \Delta)$ (see
section 3.1), we know  that when the braid $\beta(P_n; \phi)$ is
defined, it can be seen as a homotopy class of a loop in ${{\bf
D^2}\times\dots\times {\bf D^2} }\setminus \Delta$. The fundamental
group $\pi_1({{\bf D^2}\times\dots\times {\bf D^2} }\setminus \Delta)$
possesses a finite set of generators $(e_k)$. These  generators can be
expressed with the generators $\sigma_i$ of the braid group $B_n$
through the isomorphism ${\cal J}$. It follows that the proposition
will be proved once we prove that the length of the word $
\beta(P_n; \phi)$  (seen as a homotopy class written with the system of
generators $(e_k)$) is uniformely bounded when both  $P_n$ and $(\phi,
\dots ,\phi)(P_n)$ are in $\Omega^{2n}$.

\NI Consider the blowing-up set ${\cal K}$ of the generalized diagonal
$\Delta$ in $\dn$. More precisely, ${\cal K}$ is the compact set of
points:
$$
(p_1, \dots, p_i,\dots  ,p_n; p'_1,\dots, p'_j, \dots p'_n;
\Delta_{1,1}, \dots , \Delta_{i, j}, \dots , \Delta_{n,n})
$$
 where, for all $1\leq i, j\leq n$,  $\Delta_{i,j}$ is a line
 containing $p_i$ and $p'_j$. Obviously, if $p_i \neq p'_j$, then the
 line $\Delta_{i, j}$ is uniquely defined and thus ${\dn} \setminus
 \Delta$ is naturally embedded in ${\cal K}$ as an open dense set.

\NI A map $\phi$ in $\H$ yields a map $j(\phi)=(\phi,\dots,\phi)$
defined and continuous on  $\dn$. Its restriction to ${\bf
D^2}\times\dots\times {\bf D^2}$ let $\Delta$ globally invariant. We
claim that if the map $\phi$ is a $C^1$-diffeomorphism, then it can be
extended to a continuous map on ${\cal K}$. This is done as follows:
whenever for some $1\leq i, j\leq n$, we have $p_i = p'_j$,  then a
line $\Delta_{i, j}$ containing $p_i=p'_j$ is mapped to the line
$d\phi(p_i)(\Delta_{i, j})$.

\NI Furthermore, if the map $\phi$ is in $\D$ we can choose an isotopy
$\phi_t$  from the identity to $\phi$ in $\D$. Thus, the map $\Psi$:
$$\begin{array}{rrll}
\Psi:& [0, 1] \times {\cal K} & \longrightarrow & {\cal K}\\
     & (t, p_1, \dots ,  p_n) & \longmapsto & (\phi_t(p_1), \dots ,
     \phi_t(p_n))
\end{array}$$
is continuous.

\NI Let ${\cal K}^0$ be the universal cover of ${\cal K}$,  $\pi: {\cal
K}^0 \to {\cal K}$ the standard projection, and $\Psi^0 : [0, 1] \times
{\cal K} \to {\cal K}^0$ a lift of the map $\Psi$ ($\Psi= \pi
\circ\Psi^0$). The system of generators $(e_k)$ can be chosen so that
the projection $\pi$ restricted to a fundamental domain of  ${\cal
K}^0$ is a homeomorphism onto $\Omega^{2n}$. Since ${\cal K}$ is
compact $\Psi^0([0,1]\times {\cal K})$ is also compact and consequently
covers a bounded number $k(\phi, n)$ of fundamental domains of ${\cal
K}^0$. If both
  $(p_1, \dots , p_n)$ and $(\phi_1(p_1), \dots ,\phi(p_n))$ are in
  $\Omega^{2n}$, the choice of the system of generators implies that
  the length of the word $
\beta((p_1, \dots, p_n); \phi)$   written with the system of generators
$(e_k)$ is also  uniformly bounded by $k(\phi, n)$. This achieves the
proof of the proposition.
\end{pf}

\bigskip

\begin{lem} Let  $\phi$ be a map in $\D$ which  preserves the coherent
probability  measures $\lambda_1,\dots , \lambda_n$.  Then,  for $i=1$
and $i=2$ the map:
$$
\begin{array}{rll}
 {\dn}  &  \longrightarrow & \R^+\\
 P_n    &  \longmapsto     & \theta_i(\beta (P_n; \phi))
\end{array}
$$
is integrable with respect to $\lambda_1\times\dots\times \lambda_n$.

\end{lem}
\begin{pf} The integrability is a direct consequence of Proposition 2
and of the fact that $\theta_2(b) \leq (\log 3) \theta_1(b)$ for all
$b$ in $B_n$.
\end{pf}

\bigskip

\NI By applying Theorem 1 to the cocycle $\beta$ and using lemma 3, we
immediately get:

\begin{cor} Let $\phi$ be an element in $\D$ which preserves the
coherent probability  measures $\lambda_1,\dots , \lambda_n$. Then for
$i=1$ and $i=2$:

\begin{enumerate}
\item For $\lambda_1\times\dots \times \lambda_n$ a.e. $P_n$ in $\dn$
the quantity ${{1}\over {N}}\theta_i(\beta(P_n; \phi^N))$ converges
when $N$ goes to $+\infty$ to a limit $\Theta^{(i)}(P_n; \phi)$.
\item
 The map
$$
\begin{array}{rcl}
	   {\dn}  & \longrightarrow & \R^+\\
	    P_n   & \longmapsto     & \Theta^{(i)}(P_n; \phi)
\end{array}
$$
 is integrable on $\dn$ with respect to $\lambda_1\times \dots\times
 \lambda_n$.
\end {enumerate}
\end{cor}

\NI We denote  by ${\Theta}^{(i)}_{\lambda_1, \dots ,\lambda_n}(\phi)$
the integral of this function.

\begin{cor} Let $\lambda_1, \dots , \lambda_n$ and $\mu_1 , \dots ,
\mu_n$ be two sets of coherent probability measures on $\d^2$, and let
$\phi_1$ and $\phi_2$ be two elements in $\D$ which preserve the
probability measures $\lambda_1, \dots , \lambda_n$ and $\mu_1 , \dots
, \mu_n$ respectively.

\NI
Assume there exists a map $h$ in $\H$ such that :
\begin {enumerate}
\item $ h\circ \phi_1 \,=\, \phi_2\circ h$,
\item $h\ast\lambda_j\,=\, \mu_j,$ , for $j=1,\dots, n$.
\end{enumerate}
Then for $i=1$ and $i=2$:
$$
{\Theta}^{(i)}_ {\lambda_1, \dots , \lambda_n} (\phi_1) \,\,=\,\,
{\Theta}^{(i)}_{\mu_1 , \dots , \mu_n} (\phi_2).
$$
\end{cor}

\begin{pf} Consider the  cocycle $\alpha_1$:
$$
\begin{array}{rll}
	   {\dn} \times j(\HU) & \longrightarrow &  {\cal F}_n\\
	      (P_n, j(\phi_1)) & \longmapsto     & \beta(P_n; \phi_1)
\end{array}
$$

and the cocycle $\alpha_2$:
$$
\begin{array}{rll}
	   {\dn} \times j(\K) & \longrightarrow &  {\cal F}_n\\
	     (P_n, j(\phi_2)) & \longmapsto     & \beta(P_n; \phi_2)
\end{array}
$$

\NI Thanks to Theorem 2, the corollary will be proved once we show that
$\alpha_1$ and $\alpha_2$ are weakly equivalent.
Consider the homeomorphism $j(h)$ of $\dn$ into itself. It is clear
that $j(h)\ast \lambda_1\times\dots\times \lambda_n =
\mu_1\times\dots\times \mu_n$ and that:
$$
j( {\HU})= j(h)j({\K})(j(h))^{-1}.
$$
Thus it remains to prove that $\alpha_2$ and $j(h)\circ\alpha_1$ are
cohomologous.

\NI For $\lambda_1\times\dots\times\lambda_n$ a.e. $P_n$ in $\dn$, we
have:
$$
\beta(P_n; \phi_2) \,=\, \beta(P_n; h^{-1})\,\,\beta (j(h^{-1})(P_n);
\phi_1)\,\,\beta (j(\phi_1\circ h^{-1})(P_n); h).
$$
This reads :
$$
\alpha_2(P_n, j(\phi_2))\,=
\,\beta(P_n; h^{-1})\,\,h\circ \alpha_1(P_n, j(\phi_2))\,\,\beta
(j(\phi_1\circ h^{-1})(P_n)); h).
$$
 Since the map $h^{-1}$ is in $\H$ we know from Lemma 2 that the map:
$$
\begin{array}{rll}
	   {\dn}  & \longrightarrow &  {\cal F}_n\\
	      P_n & \longmapsto     & \beta(P_n; h^{-1})
\end{array}
$$
is measurable. This shows that $\alpha_2$ and $j(h)\circ \alpha_1$ are
cohomologous.
\end{pf}

\section{A discussion about these  invariants}

\subsection {The fixed points case}
In the particular case when we consider a set of fixed points
$P_n^0=(p_1^0, \dots , p_n^0)$ of $\phi\in\D$, we get, for all $N\geq
0$:
$$
\beta (P_n^0; \phi^N)\,\,=\,\,\beta (P_n^0; \phi)^N,
$$
and consequently, for $i=1$ and $i=2$, we have:
$$
\Theta^{(i)}(P_n^0; \phi)\,\,=\,\, \theta_i(\beta (P_n^0; \phi)).
$$
For the set of Dirac measures $\delta_{p^0_1}, \dots , \delta_{p^0_n}$,
which are coherent, the invariant numbers
$\Theta^{(i)}_{\delta_{p^0_1}, \dots , \delta_{p^0_n}}$ have a clear
meaning:

\begin{itemize}
\item The invariant number $\Theta^{(1)}_{\delta_{p^0_1}, \dots ,
\delta_{p^0_n}}$ is the minimum number of generators $\sigma_i$ which
are necessary to write the word
$\beta (P_n^0; \phi)$. In the particular case of a pair of fixed points
$(p^0_1, p^0_2)$, it is the absolute value of the linking number of
these two fixed points.

\item The map $\phi$ induces a map  $\phi_\star$ on the first homotopy
group of the punctered disk ${\a^2}\setminus P_n^0$. The invariant
$\Theta^{(2)}_{\delta_{p^0_1}, \dots , \delta_{p^0_n}}$  is the growth
rate of the  map $\phi_\star$. It has been shown by R. Bowen
\cite{Bowen} that this growth rate is a lower bound of the topological
entropy $h(\phi)$ of the map $\phi$ \cite{AKM}.
The assumption on the map $\phi$ to be a diffeomorphism is required in
order to get a continuous map acting on the  compact surface
$\hat{\d_n}$ obtained from ${\d_n}={\a^2}\setminus \{q_1, \dots ,q_n\}$
after blowing up the points $q_1, \dots ,q_n$ with the circles of
directions.
\end{itemize}

\subsection {The general case}
 Let  $\phi$ be a map  in $\D$. From  Proposition 2, there exists a
 positive number $K(\phi, n)$ such that for all  $P_n$ where   $\beta
 (P_n; \phi)$ is defined:
$$
{{1}\over{\log 3}}\theta^{(2)}(\beta(P_n;  \phi))
\leq\theta^{(1)}(\beta(P_n;  \phi)) \leq K(\phi, n).
$$
 Let ${\cal M}_n(\phi)$ be the set of $n$-tuples of  $\phi$-invariant,
 coherent probability measures, and let $(\lambda_1,\dots , \lambda_n)$
 be in ${\cal M}_n(\phi)$. By integrating with respect to
 $\lambda_1\times\dots\times\lambda_n$, we get:
$$
{{1}\over{\log 3}}\Theta^{(2)}_{\lambda_1,\dots , \lambda_n} (\phi)
\leq\Theta^{(1)}_{\lambda_1,\dots , \lambda_n} (\phi) \leq K(\phi, n).
$$
It follows that the quantities:
$$
\Theta^{(i)}_{ n}(\phi)\,=\, \sup_{(\lambda_1,\dots , \lambda_n)\in
{{\cal M}_n}(\phi)}\Theta^{(i)}_ {\lambda_1,\dots , \lambda_n} (\phi),
$$
are positive real numbers which are, by construction, topological
invariants. More precisely, for any pair of maps $\phi_1$ and $\phi_2$
in $\D$, which are conjugated by a map in $\H$, we have, for $i=1$,
$i=2$ and for all $n\geq 2$:
$$
\Theta^{(i)}_{ n}(\phi_1)\,=\,\Theta^{(i)}_{ n}(\phi_2).
$$
It also results from the definitions that, for $i=1$ and $i=2$ the
sequences $n\mapsto \Theta^{(i)}_{ n}(\phi)$ are non decreasing
sequences.

\NI In the next two paragraphs, we give estimates of these invariants
that generalize the fixed points case.

\subsubsection{The Calabi invariant as a minoration of the sequence
$(\Theta^{(1)}_{ n}(\phi))_n$ }
Let $\phi$ be a map in $\H$ and
consider an isotopy $\phi_t$ from the identity to $\phi$ in $\H$. To a
pair of distinct points $p_1$ and $p_2$ in $\d^2$ , we can associate a
real number  $Ang_\phi(p_1, p_2)$ which is the angular variation of the
vector $\overrightarrow{ \phi_t(p_1)\phi_t(p_2)}$ when $t $ goes from 0
to 1 (when normalizing to 1 the angular variation of a vector
 making one loop on the unit circle in the direct direction).
Since the set $\H$ is contractible, it is clear that the map $(p_1,
p_2) \mapsto Ang_\phi(p_1, p_2)$ does not depend on the choice of the
isotopy.
If $\phi $ is in $\D$, using arguments similar to the ones we used in
the proof of Proposition 2, it is easy to check (see \cite{GG}) that
the function $Ang_\phi$ is bounded where it is defined, that is to say
on ${\d^2}\times {\d^2} \setminus \Delta$.  Let $(\lambda_1,
\lambda_2)$ be a pair of $\phi$-invariant, coherent probability
measures.  The function $Ang_\phi$ is integrable on ${\a^2}\times
{\a^2} $ with respect to $\lambda_1\times \lambda_2$ and the {\it
Calabi invariant} of the map $\phi$ with respect to the pair
$(\lambda_1, \lambda_2)$ is  defined by the following integral:
$$ {\cal C}_{\lambda_1, \lambda_2}(\phi)\,=\, \int\!\int_{{\a^2}\times
{\a^2} } Ang_\phi (p_1, p_2)d\lambda_1(p_1)d\lambda_2(p_2)$$

\NI {\bf Remark:} In \cite{Calabi}, E.~Calabi defines a series of
invariant numbers associated to symplectic diffeomorphisms on
symplectic manifolds. In the particular case of area preserving maps of
the 2-disk, only one of these invariant quantities is not trivial. In
\cite{GG}, it is shown that this invariant fits with our definition in
the particular case when the invariant measures $\lambda_i$ are all
equal to the area.

\begin{prop} Let $\phi$ be a map in $\D$ which preserves a pair of
coherent probability  measures $(\lambda_1, \lambda_2)$. Then:
$$
\vert{\cal C}_{\lambda_1, \lambda_2}(\phi)\vert\,\leq \,
{\Theta}^{(1)}_{\lambda_1, \lambda_2} (\phi).
$$
Consequently, for any map $\phi$ in $\D$, and for all $n\geq 2$:
$$
\sup_{(\lambda_1, \lambda_2)\in {\cal{M}}_2(\phi)}
\vert{\cal C}_{\lambda_1, \lambda_2}(\phi)\vert
\,\leq \, {\Theta}^{(1)}_2 (\phi) \,\leq \, {\Theta}^{(1)}_n (\phi).
$$
\end{prop}

 \begin{pf}
In order to compute the Calabi invariant of $\phi$ with respect to
$(\lambda_1, \lambda_2)$, we can use the Birkhoff ergodic theorem. The
corresponding  Birkhoff sums:
$$
Ang_\phi (p_1, p_2) \,+\, Ang_\phi(\phi(p_1), \phi(p_2))\,+\,\dots +\,
Ang_\phi (\phi^{N-1}(p_1), \phi^{N-1}(p_2)),
$$
are equal to:
$$
Ang_{\phi^N}(p_1, p_2).
$$
It follows that for $\lambda_1\times \lambda_2$ almost every points
$(p_1, p_2)$ the  following limit:
$$
\tilde{Ang}_\phi(p_1, p_2) \,=\, \lim_{N\to +\infty} {{1}\over{N}}
Ang_{\phi^N}(p_1, p_2),
$$
exists and is integrable. Furthermore:
$$
\int\!\int_{{\a^2}\times {\a^2} }\tilde{ Ang}_\phi (p_1,
p_2)d\lambda(p_1)d\lambda(p_2)\,=\, \int\!\int_{{\a^2}\times {\a^2} }
Ang_\phi (p_1, p_2)d\lambda(p_1)d\lambda(p_2).$$
It is clear that for all $N\geq 0$ and for all pair of distinct points
$p$ and $q$ in $D^2$, we have the following estimate:
$$
\vert L(\beta(p, q; \phi^N))\,- \vert  Ang_{\phi^N}(p, q) \vert\,\vert
\leq 3.
$$
Dividing by $N$ and making $N$ going to $\infty$, we get for
$\lambda_1\times \lambda_2$ a.e $(p_1, p_2)$:
$$
\vert {\tilde{Ang}}_{\phi}(p_1, p_2)\vert \,\,=\,\, \Theta^{(1)}(p_1,
p_2; \phi).
$$
The integration gives:
$$
\vert{\cal C}_{\lambda_1, \lambda_2}(\phi)\vert\,\leq \,
\int\!\int_{{\a^2}\times {\a^2} }\vert\tilde{ Ang}_\phi (p_1, p_2)\vert
d\lambda_1(p_1)d\lambda_2(p_2) \,\leq \,
{\Theta}^{(1)}_{\lambda_1,  \lambda_2} (\phi).
$$
\end{pf}

\subsubsection{The topological entropy as a majoration of the sequence
$(\Theta_{2, n}(\phi))_n$ }
Let ${{\rm Diff}^\infty({\a^2}, \partial {\a^2})}$ denote the subset of
elements in $\D$ which are $C^\infty$-diffeomor\-phisms. The following
theorem can be seen as an extension  of the Bowen result (\cite{Bowen})
when we relax the hypothesis of invariance of finite sets.

\begin{thm}
Let $\phi$ be an element in ${{\rm Diff}^\infty({\d^2}, \partial
{\d^2})} $ with entropy $h(\phi)$ and $\lambda_1, \dots , \lambda_n$ a
set of coherent, $\phi$-invariant probability measures. Then, for
$\lambda_1\times\dots\times\lambda_n$-a.e.  $P_n$ in $\dn$, we have:
 $$
\Theta^{(2)}(P_n; \phi)\,\,\leq \,\, h(\phi),
$$
and consequently, for all $n\geq 2$:
$${\Theta}^{(2)}_{\lambda_1\times\dots\times\lambda_n}(\phi)
\,\,\leq\,\, h(\phi).$$
\end{thm}

\begin{pf}
Let $P_n=(p_1, \dots ,p_n)$ be a point in $\Omega^{2n}$ and
 let $h_{P_n}$ be an element in $\D$ with the following properties:

\begin{enumerate}
\item  For each $1\leq i\leq n$, $h_{P_n}$ maps $q_i$ on $p_i$.
\item  There exists an isotopy $(h_{P_n; t})_{t\in [0, 1]}$ from the
identity to $h_{P_n}$ such that, for each $1\leq i\leq n$, the arc
$\{(h_{P_n; t}, t)\,\,\,\, t\in [0, 1]\}$ coincides with the segment
joining $(q_i, 0)$ to $(p_i, 1)$ in $D^2 \times [0, 1]$.
\item The map
 $$\begin{array}{rcl}
   \Omega^{2n} & \longrightarrow & \D\\
       P_n & \longmapsto & h_{P_n}
\end{array}$$
is continuous when $\D$ is endowed with the $C^1$ topology.
\end{enumerate}

\NI Let $\phi$ be in ${{\rm Diff}^\infty({\a^2}, \partial {\a^2})} $
and assume that for $N\geq 0$, $P_n$ and $j(\phi^N)(P_n)$ are in
$\Omega^{2n}$. From the above construction , it follows that the map:
$$
\Psi^{(N)}_{P_n}\,\,=\,\, h^{-1}_{j(\phi^N)(P_n)} \circ \phi^N \circ
h_{P_n}$$ let the points $q_1, \dots , q_n$ invariant and that its
restriction to $\d_n$ induces, through the Birman isomorphism (Theorem
3),  the braid $\beta(P_n; \phi^N).$

\NI Let $P_n$ be such that the limit $\Theta^{(2)}(P_n; \phi)$ exists
(which is true for a $\lambda_1 \times \dots \times \lambda_n$ full
measure set of points). It follows that there exists $\epsilon >0$ such
that, for $N$ big enough, we have:
$$
\Theta^{(2)}(P_n; \phi) \,-\epsilon \leq{{1}\over{N}}
\theta^{(2)}(\beta(P_n; \phi^N))$$

\NI which reads:
$$(*)\,\,\, e^{(\Theta^{(2)}(P_n; \phi) \,-\epsilon)N} \leq \sup_{i=1,
\dots ,n}L(x_i\beta(P_n; \phi^N))$$

\NI It is standard that there exists a constant $c(n)>0$ such that for
any differentiable loop $\tau$ in $\hat{\d_n}$:
$$L([\tau])\,\,\leq \,\, c(n)l(\tau) $$
where $l$ stands for the euclidian length and $[-]$ is the homotopy
class in $\hat{\d_n}$. Combined with (*), this argument gives:
$$
{{1}\over{c(n)}}e^{(\Theta^{(2)}(P_n; \phi) \,-\epsilon)N} \leq
\sup_{i=1, \dots ,n}l(\Psi^{(N)}_{P_n}(x_i)).
$$
That is to say:
$$
{{1}\over{c(n)}}e^{(\Theta^{(2)}(P_n; \phi) \,-\epsilon)N} \leq
\sup_{i=1, \dots ,n}l( h^{-1}_{j(\phi^N)(P_n)} \circ \phi^N \circ
h_{P_n}(x_i))
$$
and consequently:
$$
(**)\,\,\,\, {{1}\over{c(n)}}e^{(\Theta^{(2)}(P_n; \phi) \,-\epsilon)N}
\leq\Vert h^{-1}_{j(\phi^N)(P_n)}\Vert_{1} \sup_{i=1, \dots ,n}l(
\phi^N ( h_{P_n}(x_i)),$$

\NI where $\Vert-\Vert_{1}$ stands for the $C^1$ norm.

\NI Assume  from now on that the point $P_n$ is a recurrent point of
the map $j(\phi)$ (which is true for a $\lambda_1 \times \dots \times
\lambda_n$ full measure set of points) and let $\nu(N)$ be a sequence
of integers so that:
$$
\lim_{N\to \infty}j(\phi^{\nu(N)})(P_n) \,\, =\,\, P_n.
$$
From (**) we deduce that there exists a subsequence $\hat\nu(N)$ of
$\nu(N)$ such that for $N$ big enough, and for some $1\leq i_0\leq n$,
we have:
$$
{{1}\over{c(n)K}}e^{(\Theta^{(2)}(P_n; \phi) \,-\epsilon)N}
\,\,\leq \,\, l( \phi^N ( h_{P_n}(x_{i_0})),$$
where $K$ is an upper bound of $\Vert h^{-1}_{P_n}\Vert_1$.

\NI In conclusion, the euclidian length of the loop $ h_{P_n}(x_{i_0})$
increases exponentially under iteration of the map $\phi$ with a growth
rate which is at least $\Theta^{(2)}(P_n; \phi) \,-\epsilon$. Using a
result by Y. Yomdim  (see \cite{Yomdim} and \cite{Gromov}), in the case
of $C^\infty$ maps,
we know that this provides a lower estimate of the topological entropy.
Namely, if $\phi$ is in  $ {{\rm Diff}^\infty({\a^2}, \partial
{\a^2})}$, we get:
 $$\Theta^{(2)}(P_n; \phi) \,-\epsilon\,\,\leq\,\, h(\phi).$$
Since, this is true for any $\epsilon >0$, this yields:
 $$
\Theta^{(2)}(P_n; \phi)\,\,\leq \,\, h(\phi),
$$
and by integrating:
$${\Theta}^{(2)}_ {\lambda_1\times\dots\times\lambda_n}(\phi)
\,\,\leq\,\, h(\phi).$$
\end{pf}
\vskip 2cm
\NI{\bf ACKNOWLEDGEMENTS:} It is a pleasure for the authors to thank E.
Ghys and D. Sullivan for very helpful comments and suggestions. Most
part of this work was achieved during a visit of J.M. Gambaudo at the
Institute for Mathematical Sciences. He is very grateful to the I.M.S.
for this invitation.

\vskip 2cm
\NI===================================================

\NI{\bf J.-M. Gambaudo} - {\sc Institut Non Lin\'eaire de Nice,} U.M.R.
du C.N.R.S 129,  1361 route des Lucioles, 06560 Valbonne, France.

\NI===================================================

\NI{\bf E. E.  P\'ecou} - {\sc Institute for Mathematical Sciences,}
SUNY at Stony Brook, Stony Brook, NY 11794, U.S.A.

\NI===================================================
\end{document}